\documentclass[12pt]{amsart}
\usepackage{amsmath,amscd,amssymb,amsfonts,latexsym}
\usepackage[all]{xy}
\usepackage{ulem}

\usepackage[latin1]{inputenc}

\usepackage{diagbox}

\usepackage[T1]{fontenc}

\usepackage{eucal} 
\usepackage{cmmib57} 
\usepackage{amsbsy}
\usepackage{amsthm,alltt}

\usepackage{graphicx,enumerate}

\usepackage{geometry}
\newcommand{\CP}{\mathbb{CP}^{n+1}}
\newcommand{\CN}{\mathbb{C}^{n+1}}
\newcommand{\U}{\mathcal{U}}
\newcommand{\C}{\mathbb{C}}
\newcommand{\Z}{\mathbb{Z}}

\newcommand{\Q}{\mathbb{Q}}

\newcommand{\K}{\mathcal{L}}
\newcommand{\R}{\mathcal{R}^{\bullet}}
\newcommand{\M}{\mathbb{M}}
\newcommand{\N}{\mathcal{N}}
\newcommand{\V}{\mathcal{V}}

\newcommand{\cL}{\mathcal{L}}
\newcommand{\cR}{\mathcal{R}}

\newtheorem{thm}{Theorem}[section]
\newtheorem{prop}[thm]{Proposition}
\newtheorem{lem}[thm]{Lemma}
\newtheorem{cor}[thm]{Corollary}
\theoremstyle{definition}
\newtheorem{definition}[thm]{Definition}
\newtheorem{example}[thm]{Example}
\theoremstyle{remark}
\newtheorem{remark}[thm]{Remark}
\theoremstyle{que}

\numberwithin{equation}{section}

\newcommand{\Hom}{{\rm Hom}}

\newcommand{\MHM}{{\rm MHM}}

\def\be{\begin{equation}}
\def\ee{\end{equation}}

\def\bt{\begin{thm}}
\def\et{\end{thm}}

\def\bc{\begin{cor}}
\def\ec{\end{cor}}

\def\br{\begin{remark}}
\def\er{\end{remark}}

\def\bp{\begin{prop}}
\def\ep{\end{prop}}

\def\bl{\begin{lem}}
\def\el{\end{lem}}

\def\bex{\begin{example}}
\def\eex{\end{example}}

\def\bd{\begin{definition}}
\def\ed{\end{definition}}

\title{Spectral pairs, Alexander modules, \\ and boundary manifolds}

\author{Yongqiang Liu}

\address{Y. Liu: KU Leuven, Department of Mathematics,
Celestijnenlaan 200B, B-3001 Leuven, Belgium} 

\email{liuyq1117@gmail.com}

\author{Lauren\c{t}iu Maxim}
\address{L. Maxim: Department of Mathematics, University of Wisconsin-Madison, 480 Lincoln Drive, Madison, WI 53706, USA}
\email {maxim@math.wisc.edu}

\date{\today}



\begin{document}

\begin{abstract}
Let $f: \CN \rightarrow \C $ be a reduced polynomial map, with $D=f^{-1}(0)$, $\U=\CN \setminus D$ and boundary manifold $M=\partial \U$.  Assume that $f$ is transversal at infinity and $D$ has only isolated singularities. Then the only interesting non-trivial Alexander modules of $\U$ and resp. $M$ appear in the middle degree $n$. We revisit the mixed Hodge structures on these Alexander modules and study their associated spectral pairs (or equivariant mixed Hodge numbers).
We obtain upper bounds for the spectral pairs of the $n$-th Alexander module of $\U$, which can be viewed as a Hodge-theoretic refinement of Libgober's divisibility result for the corresponding Alexander polynomials.  For the boundary manifold $M$, we show that the spectral pairs associated to the non-unipotent part of the $n$-th Alexander module of $M$ can be computed in terms of local contributions (coming from the singularities of $D$) and contributions from ``infinity''.
\end{abstract}

\maketitle


\section{Introduction}
Let $f: \CN \rightarrow \C $ be a reduced polynomial map. Set $$D=f^{-1}(0) \ \ {\rm  and} \ \  \U=\CN \setminus D.$$ The study of topology of the complement $\U$ is a classical subject going back to Zariski. The boundary manifold $M$ can be defined as the boundary of a closed regular neighborhood $\N$ of the subvariety $\CP\setminus \U$ in $\CP$, see \cite[page 149]{D1}. While the complement $\U$ has the homotopy type of a finite CW-complex of real dimension $n+1$, the boundary manifold $M$ is a smooth real closed manifold of dimension $2n+1$.

There exists a deep connection between the topology of the complement $\U$ and that of the boundary manifold $M$. For example, the embedding $M\hookrightarrow \U$ induces an $n$-homotopy equivalence, e.g., see \cite[page 150]{D1}. Moreover, Cohen-Suciu \cite{CoS1} showed that, under certain Hodge-theoretic conditions,  the cohomology ring of the complement $\U$ functorially determines that of the boundary manifold $M$.

In this paper, we focus on Hodge-theoretic aspects of the Alexander modules of $\U$ and $M$, respectively.   

In Section $2$, we review the definitions of Alexander modules  and summarize some old and new results concerning these invariants. In particular, under a certain torsion assumption (see Corollary \ref{211}), we give a new proof of a polynomial identity relating the Alexander polynomials of $\U$ and $M$.

The singularities of $D$ affect the topology of both $\U$ and $M$. Moreover, they also affect the mixed Hodge structures on the (torsion) Alexander modules of $\U$ and $M$. Assume that $f$ is transversal at infinity (see Definition \ref{ti}) and $D$ has only isolated singularities. Then the only interesting non-trivial Alexander modules of $\U$ and resp. $M$ appear in degree $n$, and there exist mixed Hodge structures on these torsion Alexander modules. In Sections $3$ and $4$, we investigate the relation between the singularities of $D$ and the corresponding  spectral pairs associated to the Alexander modules of $\U$ and $M$, respectively.  


\subsection{Alexander modules of hypersurface complements}

If $D=f^{-1}(0)$ is a plane curve \cite{L82, L83}, or a hypersurface with only isolated singularities (including at infinity) \cite{L94}, Libgober introduced and studied linking number Alexander-type invariants (i.e., induced by the polynomial $f$) associated to the hypersurface complement $\U$.  In particular, he obtained a divisibility result 
asserting that the only (possibly) non-trivial {\it global} Alexander polynomial of $\U$ divides the product of the {\it local} Alexander polynomials associated with each singular point (including at infinity). 
 
More recently, the second named author \cite{Max} used intersection homology  to provide generalizations of Libgober's results to the case of hypersurfaces with arbitrary singularities, provided that $f$ is
transversal at infinity (i.e., the hyperplane at infinity is transversal in the stratified sense to the projective completion of $D$). Furthermore, Dimca-Libgober \cite{DL} showed that for a polynomial transversal at infinity there exist canonical mixed Hodge structures on the (torsion) Alexander modules of the hypersurface complement. These mixed Hodge structures were further refined by the first named author in \cite[Theorem 1.5]{Liu} by using nearby cycles. 

\medskip

Part of this paper is devoted to the study of {\it spectral pairs} associated to the mixed Hodge structures on the Alexander modules of the complement $\U=\CN \setminus D$, provided that $f$ is transversal at infinity. 

Recall that if $A$ is a finite dimensional $\Q$-vector space endowed with a mixed Hodge structure (MHS, for short), the {\it mixed Hodge number} $h^{p,q}(A)$ is defined as the dimension of $$Gr^{p}_{F} Gr_{p+q}^{W} (A\otimes \C).$$ If, moreover, $T$ is a finite order  automorphism of the MHS $A$, we shall denote by $h_{\alpha}^{p,q}$ the dimension of the $\lambda$-eigenspace for 
the action of $T$ on $ Gr^{p}_{F} Gr_{p+q}^{W} (A\otimes \C),$ 
where $ \lambda =e^{ 2\pi i \alpha}$ and $ \alpha \in [0,1)$.  The collection $\{ h_{\alpha}^{p,q} \}$  forms the {\it spectral pairs} (or {\it equivariant mixed Hodge numbers}) of $T$ on the MHS $A$.

For example, if $F_x$ is the Milnor fiber of a hypersurface singularity germ $(D,x)$, with monodromy homeomorphism $h:F_x \to F_x$, the cohomology groups $H^i(F_x;\Q)$ carry  natural mixed Hodge structures. But the monodromy operator $h^*$ is not a MHS morphism in general. However, the semi-simple part $h^*_s$ of the monodromy operator $h^*$ is a MHS morphism of finite order, so it can be used to define the spectral pairs $h_{\alpha}^{p,q}(H^i(F_x;\Q))$ of the singularity germ $(D,x)$.

\medskip

Assume now that $f$ is transversal at infinity, and consider the infinite cyclic cover $\U^{c}$ of $\U$ associated to the epimorphism 
\begin{center}
 $f_{\ast}: \pi_{1}(\U)\twoheadrightarrow \pi_{1}(\C^{\ast})=\Z$
\end{center}  induced by $f$. Then, under the deck group action, each homology group $H_{i}(\U^{c};\Q)$ becomes a finitely generated $\Gamma:=\Q[t,t^{-1}]$-module, called the {\it $i$-th Alexander module} of the hypersurface complement $\U$. 

If $D$ has only isolated singularities, then $H_n(\U^c; \Q)$ is the only interesting non-trivial Alexander module (e.g, see Theorem \ref{Max}). It is a finite dimensional $\Q$-vector space with a canonical MHS, and it is semi-simple as a $\Gamma$-module. Moreover, the automorphism on $H_n(\U^c;\Q)$ induced by the deck transformation of $\U^c$, or equivalently the multiplication by $t$ on the $\Gamma$-module $H_n(\U^c; \Q)$, is a MHS homomorphism. So the corresponding spectral pairs are well-defined.  For simplicity, we consider the cohomology MHS on $H^{n}(\U^{c};\Q)\cong \Hom(H_{n}(\U^{c};\Q), \Q)$, where $\Q$ is regarded as a  MHS of weight $(0,0)$. 

The following theorem can be viewed as a refinement of Libgober's divisibility result in the language of spectral pairs (see Theorem \ref{t3.2}):
\bt\label{t11}   Assume that $f:\CN\to \C$ is a reduced degree $d$ polynomial transversal at infinity, so that $D=f^{-1}(0)$ has only isolated singularities. Let $\Sigma$ denote the singular locus of $D$. Then the MHS on $H^{n}(\U^{c};\Q)$ has only two possible weights: $n$ and $n+1$.  Since $H^{n}(\U^{c};\Q)$ is semi-simple under the $t$-action, the MHS on $H^{n}(\U^{c};\Q)$ splits into two pure sub-MHS, $H^{n}(\U^{c},\Q)_{\neq 1}$ and $H^{n}(\U^{c};\Q)_{= 1}$, which are pure polarized Hodge structures of weight $n$ and $n+1$, respectively. Moreover, we have the following upper bounds for the possible non-vanishing spectral pairs on $H^{n}(\U^{c};\Q)$:  
\begin{itemize}
\item[(a)] weight $n$ and $\alpha>0$, 
$$  
h_{\alpha}^{p,n-p}(H^n(\U^{c})) 
\leq  \min\lbrace \sum_{x\in \Sigma} h_{\alpha}^{p,n-p}(H^n(F_{x})); h^{p,n-p}_\alpha (H^n(F_\infty)) \rbrace,
$$
\item[(b)]  weight $n+1$ and $\alpha=0$,
$$ 
h_{\alpha}^{p,n+1-p}(H^n(\U^{c})) 
\leq  \min\lbrace \sum_{x\in \Sigma} h_{0}^{p,n+1-p}(H^n(F_{x})) + h^{p,n+1-p}(H^n(D)) ;h^{p,n+1-p}_0 (H^n(F_\infty))\rbrace. 
$$ 
\end{itemize}
where  $F_{\infty}$ is the Milnor fibre corresponding to the top degree part $f_d$ of $f$, i.e.,  $F_{\infty}=f_{d}^{-1}(1)$, and all cohomology is taken with $\Q$-coefficients.
\et 

\br The spectral pairs {\it at infinity}, $h^{p,q}_\alpha (H^n(F_\infty;\Q))$, are computed by Steenbrink \cite{Ste2} only in terms of $n$ and $d$, see Section \ref{sis}.
\er

\br It was first proved by Libgober (\cite[Proposition 3.2]{L96}) that for a hypersurface with only isolated singulariies  (including at infinity), the MHS on $H^{n}(\U^{c})$ has only weights $n$ and $n+1$.\er 
In the special case of line arrangements, the above Theorem \ref{t11} provides a refinement of a result of Massey \cite{Mas} in the language of spectral pairs (see Theorem \ref{t3.3}).


\subsection{Alexander modules of the boundary manifold}
In Section $4$, we revisit the construction of mixed Hodge structures on the Alexander modules of the boundary manifold $M$ (compare with (\cite{LiM}), and investigate the associated spectral pairs.

Consider the infinite cyclic cover $M^c$ of $M$ associated to the composition: 
$$ \pi_1(M) \twoheadrightarrow \pi_1(\U) \overset{f_*}{\twoheadrightarrow} \pi_1(\C^*)= \Z,  $$ where the first map is induced by the inclusion $M\hookrightarrow \U$ (which is an $n$-homotopy equivalence). The Alexander modules of $M$ are defined as the rational homology  groups $H_i(M^c;\Q)$, $i \in \Z$, of $M^c$, which become $\Gamma$-modules under the deck group action.


If $f$ is transversal at infinity, then it follows from Proposition \ref{p2.1} that all Alexander modules of $M$ are in fact $\Gamma$-torsion. Moreover, the zeros of the associated Alexander polynomials are roots of unity (see Theorem \ref{c2.8}), hence the semi-simple part of the $t$-action on $H_i(M^c;\Q)$, $i \in \Z$, is a finite order automorphism.
If, moreover, $D=f^{-1}(0)$ has only isolated singularities, the only interesting non-trivial Alexander module for $M$ is in the middle degree $n$. 
The second goal of this paper is  to show that  the spectral pairs associated to the non-unipotent part (i.e., corresponding to monodromy eigenvalues other than $1$) of the $n$-th Alexander module of $M$ can be computed only in terms of the {\it local} spectral pairs (at singular points of $D$) and  spectral pairs {\it at infinity}.  For simplicity,  we work again with the cohomology groups $H^{n}(M^{c};\Q)$. 
Note that the automorphism induced by multiplication by $t$ on the $\Gamma$-module $H^n(M^c;\Q)$ is not a MHS morphism in general. However, the semi-simple part of the $t$-action on $H^n(M^c;\Q)_1$ acts trivially, hence it is a MHS morphism; moreover, the following theorem asserts that the semi-simple part of the $t$-action on $H^{n}(M^{c};\Q)_{\neq 1}$ is  a MHS morphism of finite order, since  it is so for the local Milnor fibres (at the singular points of $D$) and resp. for the Milnor fiber at infinity.  Therefore, the semi-simple part of the $t$-action can be used to define the corresponding spectral pairs for $H^n(M^c; \Q)$. We then have the following result (see Theorem \ref{t4.5}):
\bt  Assume that $f:\CN\to \C$ is a reduced degree $d$ polynomial transversal at infinity so that $D=f^{-1}(0)$ has at most isolated singularities.  Let $\Sigma$ denote the singular locus of $D$. Then  we have  MHS and $\Gamma$-module isomorphisms:
$$
 H^{n}(M^{c};\Q)_{\neq 1} \cong H^{n}(F_{\infty};\Q)_{\neq 1} \oplus \big(  \bigoplus_{x\in \Sigma} H^{n}(F_{x};\Q)_{\neq 1} \big).
$$
  Therefore, for $\alpha>0$, we get:
$$ 
 h^{p,q}_{\alpha}(H^{n}(M^{c};\Q))=\sum_{x\in \Sigma} h_{\alpha}^{p,q}(H^n(F_{x};\Q)) + h^{p,q}_{\alpha}(H^n(F_{\infty};\Q))  .
$$ 
\et
We also compute the spectral pairs associated to the unipotent (i.e., eigenvalue-$1$) part of $H^n(M^c,\Q)$ in the case of hypersurfaces which are either plane curves (e.g., line arrangements) or rational homology manifolds with only isolated singularities, see  Propositions \ref{p4.7} and \ref{p4.8}, respectively.
 

\medskip

\noindent{\it Convention:} 
Unless otherwise specified, all homology and cohomology groups will be assumed to have $\Q$-coefficients.

\medskip

\textbf{Acknowledgments.} 
The authors are grateful to Julius Shaneson and Sylvain Cappell for their suggestion that one should investigate Hodge-theoretic properties of the peripheral complex, in connection to Alexander-type invariants of complex hypersurface complements.

Y. Liu was partially supported by  a FWO grant, a KU Leuven OT grant, and a Flemish Methusalem grant.  L. Maxim was partially supported by grants from NSF (DMS-1304999), NSA (H98230-14-1-0130), by a fellowship from the Max-Planck-Institut f\"ur Mathematik, Bonn, and by the Romanian Ministry of National Education, CNCS-UEFISCDI, grant PN-II-ID-PCE-2012-4-0156.


\section{Preliminaries}
In this section, we introduce the Alexander modules of a complex hypersurface complement and those of its boundary manifold, and review some old and new results concerning these invariants. 

\subsection{Alexander modules}

Let $f: \CN \rightarrow \C  $ be a reduced degree $d$ polynomial map, and set $$D=f^{-1}(0)$$ and $$\U=\CN\setminus D.$$ Since $\U$ is an affine $(n+1)$-dimensional complex variety, it has the homotopy type of a finite CW-complex of real dimension $n+1$.
 
 Let $V=\overline{D}$ be the projective completion of $D$ in $\CP$ and $H_{\infty}$ be the hyperplane at infinity so that $\CN=\CP\setminus H_{\infty}$. So $\U=\CP\setminus (V \cup H_{\infty})$.
 
 \bd\label{ti} The polynomial $f$ (or the affine hypersurface $D=f^{-1}(0)$) is called {\it transversal at infinity} if the projective closure $V$ of $D$ in $\CP$ is transversal in the stratified sense to the hyperplane at infinity $H_{\infty}$.
\ed

Consider the infinite cyclic cover $\U^{c}$ of $\U$ defined by the kernel of the {\it total linking number homomorphism}
\begin{center}
 $f_{\ast}: \pi_{1}(\U)\twoheadrightarrow \pi_{1}(\C^{\ast})=\Z$
\end{center} induced by $f$. Then, under the deck group action, each homology group $H_{i}(\U^{c})$ becomes a finitely generated $\Gamma:=\Q[t,t^{-1}]$-module, called the {\it $i$-th Alexander module} of the hypersurface complement $\U$.  Since $\Gamma$ is a principal ideal domain, to any finitely generated $\Gamma$-module one can associate an {\it order}, defined as the product of the generators associated to its torsion part, see \cite{M}. The corresponding order of $H_{i}(\U^{c})$ is called the {\it $i$-th Alexander polynomial} of $\U$, denoted by $\delta_{i}(\U,t)$.  If $H_{i}(\U^{c})$ is free or trival, then by convention we set $\delta_{i}(\U,t)=1$.

Since $\U$ has the homotopy type of a finite $(n+1)$-dimensional CW complex, $H_{i}(\U^{c})=0$ for $i>n+1$ and $H_{n+1}(\U^{c})$ is a free $\Gamma$-module. Hence the only interesting Alexander modules $H_{i}(\U^{c})$ and polynomials $\delta_{i}(\U,t)$ appear in the range $0\leq i\leq n$. 

\br \label{r2.1} In \cite[Theorem 2.10]{DN} it is shown that if the generic fiber of $f$ is connected, then the zeros of $\delta_{i}(\U,t)$ (for any $i$) are roots of unity. In particular, in this situation we have that $\delta_{n}(\U,t)=\overline{\delta_{n}(\U,t)}$, where $\overline{\delta_{n}(\U,t)}:= \delta_{n}(\U,t^{-1})$. This is indeed the case if $f$ is transversal at infinity, see \cite[Theorem 1.2]{DL}.  In fact, if $f$ is transversal at infinity, then the roots of $\delta_{i}(\U,t)$ ($i \leq n$) are roots of unity of order $d=\deg(f)$; see \cite[Theorem 4.1]{Max}. It is also shown in \cite[Section 6.2]{LiM2} that if $V \cup H_{\infty}$ is an essential hyperplane arrangement in $\CP$, then the zeros of $\delta_{i}(\U,t)$ are roots of unity if $i \leq n$. 
\er

\medskip

Let $\N$ be an open regular neighborhood of $V\cup H_{\infty}$ in $\CP$, see \cite[page 149]{D1}, and let $$\U_{0}:=\CP \setminus \N.$$ Then $\U_{0}$ is a manifold with boundary, and it is homotopy equivalent to $\U$. 
\bd The {\it boundary manifold of $\U$} is the $(2n+1)$-dimensional real closed manifold 
$$M:=\partial \U_{0}.$$
\ed
The inclusion $M \hookrightarrow \U$ is an $n$-homotopy equivalence (cf. \cite[(5.2.31)]{D1}). So we have an epimorphism:  \begin{center}
$ \pi_{1}(M)\twoheadrightarrow \pi_{1}(\U)\overset{f_*}{\twoheadrightarrow}  \pi_{1}(\C^{\ast})=\Z$,
\end{center}   which defines the (linking number) infinite cyclic cover $M^{c}$ of $M$. The homology group $H_{i}(M^{c})$  becomes a finitely generated $\Gamma$-module, called the {\it $i$-th Alexander module} of the boundary manifold $M$.  The Alexander polynomials $\delta_i(M,t)$  of $M$ are defined in the same way as the ones for $\U$.

\br Since the homomorphisms used to define the infinity cyclic covers for $\U$ and resp. $M$ are surjective, $\U^{c}$ and $M^{c}$ are connected topological spaces, hence $$H_{0}(\U^c) \cong \Gamma/(t-1) \cong H_{0}(M^{c}).$$
\er

By analogy with the local situation, in the rest of the paper we call the $\Gamma$-action $t: H_{i}(\U^{c})\rightarrow H_{i}(\U^{c})$ corresponding to the generator of the deck group $\Z$ of $\U^c$ the {\it monodromy action} on the Alexander module $H_{i}(\U^{c})$, and similarly for the boundary manifold $M$.

\subsection{Peripheral complex}
 Consider the local system $\K $ on the hypersurface complement $\U$ with stalk $\Gamma=\Q[t,t^{-1}]$ and representation of the fundamental group defined by the composition: 
$$ \pi_{1}(\U) \overset{f_*}{\rightarrow} \pi_{1}(\C^{\ast})   \rightarrow   Aut(\Gamma),$$
with the second map given by  $1_{\Z}\mapsto t$. 
Here we denote by $t$ the automorphism of $\Gamma$ given by multiplication by $t$. $\K$ shall be referred to as the {\it (total) linking number local system} on $\U$. Via the embedding $M \hookrightarrow \U$, one can of course restrict the local system $\K$ to $M$.

The relation between the linking number local system $\K$ and the Alexander modules of $\U$ and $M$  is given as follows: for any $i$, there exist  $\Gamma$-module isomorphisms: 
 \be H_{i}(\U^{c}) \cong H_{i}(\U,\K), \ee
 \be H_{i}(M^{c}) \cong H_{i}(M,\K|_M). \ee

\medskip 
 
 For any complex algebraic variety $X$ and any Noetherian commutative ring $R$, we denote by  $D^{b}_{c}(X,R)$ the derived category of bounded cohomologically constructible complexes of sheaves of $R$-modules on $X$. If $A^{\bullet} \in D^{b}_{c}(X, \Gamma)$, let $\mathcal{D} A^{\bullet}$ denote its Verdier dual and let $\overline{A^{\bullet}}$ be the complex obtained from $A^{\bullet}$ by composing all $\Gamma$-module structures with the involution $t\rightarrow t^{-1}$. Then we have that:
\be\label{d} \mathcal{D} \K  \simeq \overline{\K}[2n+2].\ee

\bd\label{defpc} Let $j:\U \hookrightarrow \CP$ denote the inclusion map. The {\it peripheral complex} $\R \in D_{c}^{b}(\CP, \Gamma)$ of $f$ is defined by the distinguished triangle 
  \begin{equation}
  j_{!}\K \rightarrow Rj_{\ast}\K \rightarrow \R  \overset{[1]}{\rightarrow}.
  \end{equation}
  \ed
It is clear that $\R$ has compact support on $V\cup H_{\infty}$, and 
\be\label{sup} \R \cong  (Rj_{\ast} \K)\vert_{V\cup H_{\infty}}.\ee 
Furthermore, it follows from (\ref{d}) that, up to a shift,  $\R$ is a self-dual complex, in the sense that $$ \mathcal{D}\R \simeq \overline{\R}[2n+1] .$$

\br  If $f$ is transversal at infinity, then it follows by results of the second author \cite{Max} that the peripheral complex $\R$ as defined here corresponds to the shifted Cappell-Shaneson peripheral complex $\R[-2n-2]$, which was defined as the cone of a map of certain intersection cohomology complexes, see \cite{CS} or \cite{Max} for more details. 
\er

The peripheral complex $\R$ can be regarded as a generalization of Deligne's nearby cycle complex $\psi_{f}{\Q_{\CN}}$ associated to the polynomial $f$.  Indeed, let 
\be\label{for}{\rm\it for}:  D^{b}_{c}(D,\Gamma)\rightarrow D^{b}_{c}(D,\Q)\ee
be the forgetful functor, 
which sends a {\it torsion} $\Gamma$-module sheaf complex to its underlying $\Q$-complex.  It can be shown that $\R\vert_{D}$ is always a torsion $\Gamma$-module sheaf complex (e.g., see \cite{Liu, Max}), and the following non-canonical  quasi-isomorphism (obtained by combining \cite[page 13]{Br}, \cite[Lemma 3.4(b)]{Bu} and \cite[Lemma 2.9]{LiM}) holds in $D^{b}_{c}(D,\Q)$:
\be\label{pn} {\rm\it for} (\R\vert_{D})\simeq \psi_{f}{\Q_{\CN}}[-1].\ee 
For simplicity, in the following we write $\Q$ for the constant sheaf $\Q_{\CN}$ on $\CN$.
 
\medskip
 
The nearby cycle complex $\psi_{f}{\Q}$ has many good properties. For example, $\psi_{f}{\Q}[-1]$ is a  perverse sheaf of $\Q$-vector spaces on $D$. Moreover, Deligne's nearby cycle functor $$\psi_f[-1]:D^b_c(\CN,\Q) \to D^b_c(D,\Q)$$ preserves perverse $\Q$-sheaves and can be lifted to the category of Saito's mixed Hodge modules. Due to its relation to the nearby cycles, it would be desirable for the peripheral complex $\R$  to exhibit similar properties, assuming it is a torsion $\Gamma$-module sheaf complex. However, it is known that  
such a torsion assumption does not hold in general for the peripheral complex. Nevertheless, we can single out two interesting situations when we have a positive answer:
\begin{itemize} 
\item[(a)]  $f$ is transversal at infinity;
\item[(b)]  $V\cup H_{\infty}$ defines an essential hyperplane arrangement in $\CP$.
\end{itemize}
The situation described in (a) was considered in detail in \cite[Section 3]{Max}, while the relevant statements in case (b) can be derived from the proof of \cite[Proposition 6.8]{LiM2}. The remaining of this paper will focus mainly on the case (a). 

\medskip

For future reference, if $\R$ is a torsion sheaf complex of $\Gamma$-modules,  let us set 
\be\label{notr} R^{\bullet}:={\rm\it for} (\R) [1].\ee  
Then  $R^{\bullet}[n]$ is a self-dual $\Q$-complex with support on $V\cup H_{\infty}$, i.e.,  \be\label{dpr} \mathcal{D}( R^{\bullet}[n])\simeq R^{\bullet}[n]. \ee

\medskip

The peripheral  complex has a very nice geometric interpretation, in the sense that its (hyper)cohomology realizes the Alexander modules of the boundary manifold $M$. More precisely, we have the following result.
\bp \label{p2.1} For any $i$, there are $\Gamma$-module isomorphisms:  \begin{center}
$H_{i}( M^{c})\cong H^{2n+1-i}(V\cup H_{\infty};\R)$.
\end{center} 
\ep
\begin{proof}
 The claim follows from the fact that the homology exact sequence for the pair $(\U^{c}, M^{c})$ can be identified (by Poincar\'e duality and homotopy equivalence) with the hypercohomology long exact sequence associated to the distinguished triangle defining the peripheral complex $\R$. A detailed proof for the case when $f$ is transversal at infinity can be found in \cite[Proposition 6.1]{LiM}, and  the arguments in loc.cit. can be easily adapted to the general case.
\end{proof}

The following result is a direct consequence of the local structure theorem for hypersurface singularity germ complements, cf. works of Libgober \cite{L09} and Budur-Wang \cite{BW}.
\bt\label{c2.8} If the peripheral complex $\R$ is $\Gamma$-torsion, then the zeros of the Alexander polynomials $\delta_{i}(M,t)$ (for any $i$) of the boundary manifold $M$ are roots of unity.
\et
\begin{proof}
By Proposition \ref{p2.1}, the peripheral complex $\R$ realizes the Alexander modules of the boundary manifold $M$. If $\R$ is a torsion $\Gamma$-complex, it follows from the hypercohomology spectral sequence that one only needs to check that the zeros of orders of  the stalk cohomology groups of $\R$ at points in $V \cup H_{\infty}$ are roots of unity.
Indeed, by using the compactly supported hypercohomology long exact sequence for the inclusion of strata of $V\cup H_{\infty}$, one can first reduce the problem to showing that for any stratum $S$ of $V\cup H_{\infty}$, the zeros of the order of $H^{\ast}_{c}(S,\R\vert _{S})$ are roots of unity. Then a hypercohomology spectral sequence argument reduces the problem to studying the orders of the stalk cohomologies $\mathcal{H}^i(\R)_{x\in S}$ of $\R$ at points in the stratum $S$.

If $x\in D$, we have by (\ref{pn})  the $\Gamma$-module isomorphisms
$$\mathcal{H}^i(\R)_x \cong \mathcal{H}^{i-1}(\psi_{f} \Q)_x\cong H^{i-1}(F_x),$$
for $F_x$ the local Milnor fiber at $x$. 
Hence the claim follows in this case from the classical monodromy theorem. If $x\in H_{\infty} \setminus V$, then $\mathcal{H}^i(\R)_x$ is isomorphic as a $\Gamma$-module to either $\Gamma/(t^{d}-1)$ or $0$, where $d$ is the degree of the reduced polynomial $f$. So the only case left for consideration is when $ x\in V\cap H_{\infty}$. 

If $x\in V\cap H_{\infty}$, let $\U_{x}$ denote the local complement of $x$ in $\U$, i.e., $\U_x=\U \cap B_x$, for $B_x$ a small Milnor ball centered at $x$ in $\CP$.  Note that $$\mathcal{H}^i(\R)_x\cong H^{i}(\U_{x}, Rj_{\ast} \K) \cong \overline{H_{i-1}(\U_{x}^{c})},$$ where $\U_{x}^{c}$ denotes the induced infinite cyclic cover of $\U_{x}$, and the second isomorphism follows from the Universal Coefficient Theorem for the principal ideal domain $\Gamma$ (e.g., see see \cite[2-5]{LiM}). Here $\overline{\ast}$ denotes the composition of $\Gamma$-module structures with the involution $t\rightarrow t^{-1}$.

Next we need to recall the definition of {\it jumping loci}. Let $X$ be a connected finite CW-complex with $\pi_{1}(X)=G$. Then the group of $\C$-valued  characters, $\Hom(G,\C^{\ast})$, is a commutative, affine algebraic group. Each character $\rho \in \Hom(G,\C^{\ast})$ defines a rank-one local system on $X$, denoted by $\K_{\rho}$.
The homology jumping loci of $X$  are the subsets of the character group $\Hom(G,\C^{\ast})$ defined as:
$$\V^i_{k}(X)=\lbrace \rho\in \Hom(G,\C^{\ast}) \mid \dim_{\C} H_{i}(X, \K_{\rho})\geq k \rbrace.$$ 

 The local structure theorem of  Libgober \cite{L09} and Budur-Wang \cite{BW} asserts that each $\V^{i}_{k}(\U_{x})$ is a finite union of torsion translated subtori.
As noted in \cite[Section 6.2, Proposition 6.6]{LiM2}, the zeros of the order of $H_{i}(\U_{x}^{c})$ (i.e., the roots of the local Alexander polynomial at $x$) can be obtained from $\V^{i}_{1}(\U_{x})$ by intersecting it with the $1$-dimensional torus defined by $(t^{-d}, t,\cdots,t)$. The torsion assumption for $\R$ implies that this intersection is either a finite set of points or it is the empty set. The structure theorem then yields that these points are all torsion points, hence the corresponding zeros are roots of unity.   
\end{proof}

\br It can be seen from the above proof that the condition of $\R$ being $\Gamma$-torsion is equivalent to the more geometric condition which requires the $\Gamma$-torsion property for the {\it local} Alexander modules $H_k(\U^c_x)$ at points $x \in V \cap H_{\infty}$, for all $k \in \Z$. For example, this condition is satisfied if $f$ is transversal at infinity (see \cite{Max}) or if $V \cup H$ is an essential hyperplane arrangement (see \cite{LiM2}), so the statement of Theorem \ref{c2.8} is in particular applicable to these situations.
\er


\subsection{An Alexander polynomial identity}

Let us now assume that $\R$ is a $\Gamma$-torsion complex. Then Proposition \ref{p2.1} yields that $H_{i}(M^{c})$ is a torsion $\Gamma$-module for any $i$. The $n$-homotopy equivalence induced by the embedding $M\hookrightarrow \U$ implies that there exist isomorphisms of $\Gamma$-modules \be\label{ison}
 H_{i}(M^{c}) \cong H_{i}(\U^{c}), \  \ {\rm for} \  \ i \leq n-1, 
\ee 
and a surjective $\Gamma$-module homomorphism $$H_{n}(M^{c}) \twoheadrightarrow  H_{n}(\U^{c}).$$ 
Hence we reprove the following well-known fact (e.g., see \cite{Max,Liu,LiM,LiM2}):
\bc\label{210} If the peripheral complex $\R$ is a $\Gamma$-torsion complex, then the Alexander module
$H_{i}(\U^{c})$  is a torsion $\Gamma$-module for any $i \leq n$. Moreover, $H_{n+1}(\U^{c}) \cong \Gamma^{\mu}$, where $\mu=\vert \chi(\U)\vert$.
\ec 
It is worth mentioning here that in both situations (a) and (b) mentioned before, the affine hypersurface $D$ is homotopy equivalent to a bouquet of exactly $\mu$ $n$-spheres, i.e., $$D\simeq \bigvee_\mu S^{n},$$ see \cite[page 476]{DP}, \cite[Corollary 2.2]{D3} and \cite[Proposition 2.1]{DJL}.

\medskip

Poincar\'e duality and the Universal Coefficients Theorem for the principal ideal domain $\Gamma $ yield that $$H_{n+1}(\U^{c}, M^{c}) \cong \Gamma^{\mu} \oplus \overline{H_{n}(\U^{c})}.$$ Then by using the fact that $H_{n+1}(M^c)$ is a torsion $\Gamma$-module, together with the isomorphisms (\ref{ison}), we obtain a long exact sequence of $\Gamma$-modules:
\be  \label{2.6}
0 \to \Gamma^{\mu} \overset{\rho}{\to} \Gamma^{\mu}\oplus \overline{ H_{n}(\U^{c})} \to H_{n}(M^{c}) \to H_{n}(\U^{c})\to 0,
\ee
where $\rho \in {\rm GL}_\mu(\Gamma)$ is an invertible $\mu \times \mu$ matrix. 
Let $$e(t):= \det(\rho) \in 
\Gamma,$$ which shall be referred to as the {\it error term}; see Remark \ref{r4.4} for a justification of this terminology.
We then have the following Alexander polynomial identity:
\bc  \label{211} If $\R$ is $\Gamma$-torsion, then 
\be\label{Al} 
 e(t) \cdot \delta_{n}(\U, t) \cdot \overline{\delta_{n}(\U,t)} = \delta_{n}(M, t).
\ee 
\begin{proof}
Formula (\ref{Al}) is a direct consequence of the long exact sequence (\ref{2.6}).  In fact, (\ref{2.6}) yields a short exact sequence of $\Gamma$-modules:
$$ 0 \to  \mathrm{coker}( \rho) \oplus \overline{ H_{n}(\U^{c})} \to H_{n}(M^{c}) \to H_{n}(\U^{c})\to 0.$$ Then $\mathrm{coker}( \rho)$ (possibly 0) is a torsion $\Gamma$-module. Since $\Gamma$ is a principal ideal domain, one can consider the corresponding Smith normal form for the matrix $\rho$. It is clear that the order of $\mathrm{coker}( \rho)$ coincides with the determinant of $\rho$, up to multiplication by units. Then the claim follows from the additivity property for  the order of torsion $\Gamma$-modules  associated to the short exact sequence.
\end{proof}
\ec
Note that if we are in the situations (a) or (b), then $\delta_{n}(\U,t)=\overline{\delta_{n}(\U,t)}$, so in these cases formula (\ref{Al}) simplifies to: 
\be 
 e(t) \cdot \delta_{n}(\U, t)^{2} = \delta_{n}(M, t).
\ee 
\br  This Alexander polynomial identity can also be proved by using the Reidemeister torsion, e.g., see \cite{CF,LiM}.
\er


\subsection{Self-duality and immediate consequences}

If $\R$ is $\Gamma$-torsion, then Proposition \ref{p2.1} and the self-duality of $\R$  
imply that $M^{c}$ has the homology properties of a closed real manifold of dimension $2n$.  More precisely, we recover the Duality Theorem of Milnor \cite[Section 4]{M} in the special case of the boundary manifold $M$. We recall here Milnor's result:

\bt(Duality Theorem \cite{M}) Assume that $M$ is a closed manifold of real dimension $(2n+1)$. If $H^{i}(M^{c})$ is $\Gamma$-torsion for any $i$, then $H^{2n}(M^{c})$ is one-dimensional over $\Q$, and the vector spaces $H^{i}(M^{c})$ and $H^{2n-i}(M^{c})$ are dual to each other, being orthogonally paired to $H^{2n}(M^{c})=\Q$ by the cup product pairing.
\et 

By the Duality Theorem, the intersection pairing for the middle degree term $H^{n}(M^{c})$ is non-degenerated. If $n$ is odd, then the intersection form is skew-symmetric and it is determined by 
the dimension of $H^{n}(M^{c})$. In particular, $\dim H^{n}(M^{c})$ is always even. On the other hand, if $n$ is even,  then the intersection form is symmetric, and we can diagonalize it.  Let $e_{+}$ and $e_{-}$ denote the  number of positive and resp. negative diagonal entries of the diagonalization. Then $\dim H^{n}(M^{c})=e_+ +e_-$. The {\it signature} of $M^c$ is defined as $$\sigma(M^c):=e_+ -e_-.$$

\br A natural question to ask when $n$ is even is the following: if $V \cup H_{\infty}$ is an essential hyperplane arrangement, is the signature $\sigma(M^c)$ determined by the intersection lattice of the arrangement?   
\er 

We end this section with an application of the Duality Theorem.  First, recall the following result from \cite{Max} (where the second part of the statement is based on \cite{L94}).

\bt \label{Max} Assume that the reduced degree $d$ polynomial $f: \CN\rightarrow \C$ is transversal at infinity. Then, for any $i\leq n$, the zeros of the Alexander polynomial $\delta_i(\U,t)$ associated to $H_{i}(\U^{c})$ are roots of unity of order $d$, and $H_{i}(\U^{c})$ is a semi-simple $\Gamma$-module.  Moreover,   if $V=\overline{D}$ has no codimension-one singularities, then $H_{i}(\U^{c})=0$ for $0<i<n-k$, where $k$ is the dimension of the singular locus of $D$ and (by convention) $k=-1$ if $D$ is smooth. 
\et

We can now prove the following:
\bt \label{LiM} Assume that the reduced degree $d$ polynomial $f: \CN\rightarrow \C$ is transversal at infinity. Then, for any $i\neq n$, the zeros of the Alexander polynomial  $\delta_i(M,t)$ associated to $H_{i}(M^{c})$ are roots of unity of order $d$, and $H_{i}(M^{c})$ is a semi-simple $\Gamma$-module.  Moreover,  if $V=\overline{D}$ has no codimension-one singularities, then $H_{i}(M^{c})=0$ for $0<i<n-k$ or $n+k<i<2n$, where $k\geq 0$ is the dimension of the singular locus of $D$. 
\et
\begin{proof}
As shown in \cite{Max}, since $f$ is transversal at infinity we have that $\R$ is a $\Gamma$-torsion complex. 
Next recall that the inclusion $M\hookrightarrow \U$ is an $n$-homotopy equivalence, thus \begin{center}
$H_{i}(M^{c})\cong H_{i}(\U^{c})$ for $0\leq i\leq n-1$.
\end{center} 
Then the claim for $i\leq n-1$ follows from the previous theorem. The other half of the story (including the desired vanishing) follows from the self-duality of $\R$ (or by Milnor duality).  
\end{proof}
\br  Even when the affine hypersurface $D$ is smooth, the middle Alexander module $H_{n}(M^{c})$ of the boundary manifold is never zero (unlike $H_{n}(\U^{c})$, which vanishes in this case), see Proposition \ref{p4.3} below. The collection of zeros of the corresponding Alexander polynomial $\delta _{n}(M,t)$ is made of local contributions, as well as contributions coming from ``infinity'', compare also with \cite[Theorem 1.2]{LiM}.
\er

\section{Spectral pairs of Alexander modules of hypersurface complements}

From now on, we will assume that $f$ is transversal at infinity, i.e., $V$ intersects $H_{\infty}$ transversally (in the stratified sense).

\subsection{Spectral pairs}
\bd A mixed Hodge structure (MHS, for short) is a finite dimensional $\Q$-vector space $A$ endowed with a finite increasing weight filtration $W_{\bullet}$, and with a finite decreasing Hodge filtration $F^{\bullet}$ on $A\otimes \C$ such that $( Gr_{k}^{W} A, F^{\bullet})$ is a pure Hodge structure of weight $k$, for all $k$. The {\it mixed Hodge number} $h^{p,q}(A)$ is defined as the dimension of $ Gr^{p}_{F} Gr_{p+q}^{W} (A\otimes \C).$  A morphism of mixed Hodge structures is a $\Q$-linear map between two mixed Hodge structures, which is compatible with both filtrations $F^{\bullet}$ and $W_{\bullet}$.
\ed

\bd
Assume that $T$ is a finite order automorphism of the MHS $A$, say $T^m=Id_A$. We shall denote by $h_{\lambda}^{p,q}$ the dimension of the $\lambda$-eigenspace for 
the action of $T$ on $$ Gr^{p}_{F} Gr_{p+q}^{W} (A\otimes \C).$$ 
 It is clear that $h_{\lambda}^{p,q}=0$ if $\lambda^{m}\neq 1$.
Set $ \lambda =e^{ 2\pi i \alpha}$, where $ \alpha \in [0,1)$.  In what follows, $h_{\lambda}^{p,q}$ will also  be denoted by $h_{\alpha}^{p,q}$. The collection $\{ h_{\alpha}^{p,q} \}$  forms the {\it spectral pairs} (or {\it equivariant mixed Hodge numbers}) of the MHS $(A,T)$.
\ed

\bd We say that a triple $(A, T, n)$  (with $(A,T)$ as in the above definition) is {\it symmetric} if there there is a nilpotent MHS morphism $N: A \to A$ of type $(-1,-1)$ and commuting with $T$ such that the weight filtration of $A$ is the weight filtration associated to the nilpotent operator $N$ with center $n$. 
\ed 

\br (\cite[Remark 2.8]{D4}) \\
(1) By conjugation, we get that $h^{p,q}_{\lambda} =h^{q,p}_{\overline{\lambda}}$.\\
(2) One can read the size of the Jordan blocks for $N$ from the spectral pairs in the case of a symmetric MHS. \\
(3) For a short exact sequence of MHS with compatible automorphisms, if the semi-simple part of these 3 automorphisms are all MHS morphisms of finite order, then the corresponding spectral pairs are additive. 
\er

\bex Let $F_x$ be the Milnor fiber of a hypersurface singularity germ $(D,x)$, with monodromy homeomorphism $h:F_x \to F_x$.
The cohomology group $H^i(F_x)$ carries a natural mixed Hodge structure, but the monodromy operator $h^*$ is not a MHS morphism in general. However, the semi-simple part $h^*_s$ of the monodromy operator $h^*$ is a MHS morphism of finite order, so it can be used to define the spectral pairs $h_{\alpha}^{p,q}(H^i(F_x))$ of the singularity germ $(D,x)$. Let $N$ be the logarithm of the unipotent part of the monodromy action $h^*$ on $H^{i}(F_{x})$. Then $N$ is a MHS morphism of type $(-1,-1)$.
If $f$ defines an isolated hypersurface singularity at $x$, then the triples $(H^{n}(F_{x})_{\neq 1}, h^*_s, n)$ and $(H^{n}(F_x)_{ 1}, h^*_s, n+1)$ are two symmetric MHS.  Here $H^{n}(F_x)_{1}$ and $H^{n}(F_x)_{\neq 1}$ denote the monodromy eigenspace for the eigenvalue $1$ (i.e., the unipotent part) and, resp.,  for eigenvalues other than $1$ (i.e., the non-unipotent part). 
\eex

Another example of a mixed Hodge structure with a finite order automorphism is provided by the Alexander modules of hypersurface complements, with defining polynomial transversal at infinity. Indeed, it was shown in \cite{DL} that in this case there exist canonical mixed Hodge structures on the Alexander modules $H_{i}(\U^{c})$ of the hypersurface complement, for all $i\leq n$. These mixed Hodge structures were further refined by the first author in \cite{Liu} by using nearby cycles.  
Let us summarize here the relevant results from \cite{DL} and \cite{Liu}. For simplicity, we consider the cohomology MHS on $H^{i}(\U^{c})= \Hom(H_{i}(\U^{c}), \Q)$, where $\Q$ is regarded as a  MHS of weight $(0,0)$. 

\bt  \label{t3.1}
 Assume that the reduced degree $d$ polynomial $f: \CN\rightarrow \C$ is transversal at infinity. Then there exists a canonical MHS on $H^{i}(\U^{c})$  ($i\leq n$), which is compatible with the action $t$ of $\Gamma$, i.e., the monodromy $t: H^{i}(\U^{c})\rightarrow H^{i}(\U^{c})$ is a MHS morphism.  Moreover, there exist MHS isomorphisms \be \label{3.0}
H^{i}(F_{\infty}) \cong H^{i}(D, \psi_{f}\Q) \cong H^{i}(\U^{c}) \text{ for } i<n,
\ee and two injective MHS morphisms: \be \label{3.1}
 H^{n}(\U^{c}) \hookrightarrow H^{n}(F_{\infty}) 
 \ee  and 
 \be  \label{3.2}
 H^{n}(\U^{c})\hookrightarrow H^{n}(D, \psi_{f}\Q),
\ee compatible with the respective monodromy actions, where $F_{\infty}$ is the Milnor fibre of the top degree part $f_d$ of $f$, i.e.,  $F_{\infty}=\lbrace f_{d}=1 \rbrace$. 
\et

In particular, we can define spectral pairs associated to the Alexander modules $H^{i}(\U^{c})$  ($i\leq n$) of hypersurface complements.


\subsection{Isolated singularities}\label{sis} In this section, we study in detail the spectral pairs of the Alexander modules of complements to hypersurfaces with only isolated singularities.

Let $f:\CN \to \C$ be a polynomial transversal at infinity, so that the hypersurface $D:=f^{-1}(0)$ has only isolated singularities. Let $\Sigma$ denote the singular set of $D$. 
By Theorem \ref{Max}, the only interesting (possibly non-trivial) Alexander module of the hypersurface complement $\U=\CN\setminus D$ is $H_{n}(\U^{c})$. Moreover, as shown by Libgober in \cite{L94} (see also \cite{Max} and \cite{Liu}), the {\it global} Alexander polynomial $\delta_n(\U,t)$ can be described in terms of similar {\it local} Alexander polynomials associated to the local link complements at each singular point. In fact, Libgober's divisibility result can be derived from the two injective maps appearing in Theorem \ref{t3.1}. More precisely, (\ref{3.1}) yields that (e.g., see the proof of \cite[Corollary 7.5]{LiM}) 
\be\label{d1}
\delta_{n}(\U,t)\mid (t-1)^{(-1)^{n+1}}(t^{d}-1)^{\xi}
\ee 
where $\xi=\dfrac{(d-1)^{n+1}+(-1)^{n}}{d}$, while (\ref{3.2}) leads to (see \cite[Section 5.2.1]{Liu})
\be\label{d2}
\delta_{n}(\U,t)\mid (t-1)^{\mu}\prod_{x\in \Sigma}\Delta_{x}(t),
\ee 
where $ \Delta_{x}(t)$ is the (top) local Alexander polynomial associated to the singular point $x\in \Sigma$. Moreover, in this case we have: $\mu=(d-1)^{n+1} - \sum_{x\in \Sigma} \mu_{x}$, where $\mu_{x}$ is the Milnor number of $f$ at $x$.

\medskip

The purpose of this section is to refine the divisibility results (\ref{d1}) and (\ref{d2}) in the language of spectral pairs. Note that whenever we focus on the isolated singularity case,  the only (possibly) non-trivial cohomology group for any of $\U^{c}$, $M^{c}$, the local Milnor fibre $F_x$ at a singular point $x \in \Sigma$, the Milnor fibre $F_{\infty}$ at infinity, and even for the affine hypersurface $D$, appears in degree $n$. As a notational convention, if $X$ is any of these spaces, in what follows we simply write $h_{\alpha}^{p,q}(X)$ for $h_{\alpha}^{p,q}(H^{n}(X))$.

 Since $D$ is assumed to have only isolated singularities, the transversality assumption for $f$ implies that $f_{d}$ is a degree $d$ homogeneous polynomial which defines an isolated singularity at the origin in $\CN$. Let $F_\infty=\{f_{d}=1\}$ be the Milnor fiber of $f_{d}$ as in Theorem \ref{t3.1}. Then $H^n(F_\infty)$ has a MHS so that the corresponding monodromy operator $h^*:H^n(F_\infty) \to H^n(F_\infty)$ is a MHS morphism (since $h$ is an algebraic map in this case). In fact, Steenbrink \cite{Ste2} showed that the MHS on $H^{n}(F_{\infty})$ has only weights $n$ and $n+1$, and he computed the spectral pairs $$h_{\alpha}^{p,q}(F_{\infty}):=h_{\alpha}^{p,q}(H^n(F_{\infty}))$$ as we shall now recall. 
 
Since the spectral pairs do not depend on the choice of $f_{d}$ (i.e., they are computable only in terms on $n$ and $d$), one can just take $f_{d}=x_{0}^{d}+\cdots +x_n^d .$
Let $\M(n,d)$ denote the Milnor algebra associated to the homogeneous polynomial $x_0^d+\cdots +x_n^d$. Then $\M(n,d)$ has a finite rank as a $\C$-vector space, and we can take a system of monomials as a basis for $\M(n,d)$. We denote by $\M(n,d)_{m}$ the subalgebra of $\M(n,d)$ generated by the monomials of degree $m$ (modulo the ideal $(x_{0}^{d-1}, \cdots, x_{n}^{d-1})$) in the basis. Then the following formulae hold:
\be 
h_{\alpha}^{p,n-p}(F_{\infty}) =\left\{ \begin{array}{ll}
0, & \alpha=0, \\
\dim \M(n,d)_{pd-n-1+d\alpha}, & \alpha>0.\\
\end{array}\right.
\ee
\be\label{36}
h_{\alpha}^{p,n+1-p}(F_{\infty}) =\left\{ \begin{array}{ll}
\dim \M(n,d)_{pd-n-1}, & \alpha=0, \\
0, & \alpha>0.\\
\end{array}\right.
\ee

The following result is a consequence of Theorem \ref{t3.1} and Steenbrink's computation. 

\bt \label{t3.2}  Assume that $f:\CN\to \C$ is a reduced degree $d$ polynomial transversal at infinity, so that $D=f^{-1}(0)$ has only isolated singularities. Let $\Sigma$ denote the singular locus of $D$. Then the MHS on $H^{n}(\U^{c})$ has only two possible weights: $n$ and $n+1$.  Since $H^{n}(\U^{c})$ is semi-simple, the MHS on $H^{n}(\U^{c})$ splits into two pure sub-MHS, $H^{n}(\U^{c})_{\neq 1}$ and $H^{n}(\U^{c})_{= 1}$, which are pure polarized Hodge structures of weight $n$ and $n+1$, respectively. Moreover, we have the following upper bounds for the possible non-vanishing spectral pairs of $H^{n}(\U^{c})$:  
\begin{itemize}
\item[(a)] weight $n$ and $\alpha>0$, \be  \label{3.3}
h_{\alpha}^{p,n-p}(\U^{c}) 
\leq  \min\lbrace \sum_{x\in \Sigma} h_{\alpha}^{p,n-p}(F_{x}); \dim \M(n,d)_{pd-n-1+d\alpha} \rbrace,
\ee
\item[(b)]  weight $n+1$ and $\alpha=0$,
\be \label{3.4}
h_{\alpha}^{p,n+1-p}(\U^{c}) 
\leq  \min\lbrace \sum_{x\in \Sigma} h_{0}^{p,n+1-p}(F_{x}) + h^{p,n+1-p}(D) ;\dim \M(n,d)_{pd-n-1}\rbrace. 
\ee 
\end{itemize}
\et 
\begin{proof}

The claim about possible weights on $H^n(\U^c)$ follows from  (\ref{3.1}) and Steenbrink's computation. 
In view of  (\ref{3.2}), we also need to investigate the MHS on $H^{n}(D,\psi_{f}\Q)$. 
 
We have a distinguished triangle in the bounded derived category of mixed Hodge modules on $D$: \be\label{dist1} \Q_{D} \to \psi_{f}\Q \to \varphi_{f}\Q \overset{[1]}{\to}.\ee There are decompositions 
$\psi_{f}\Q=\psi_{f, 1}\Q \oplus \psi_{f,\neq 1}\Q$, and similar for the vanishing cycles $\varphi_{f}\Q$, so that $\psi_{f,\neq 1}\Q\simeq \varphi_{f,\neq 1}\Q$ and the semi-simple part of the corresponding monodromy operators act trivially on $\psi_{f, 1}\Q$ and $\varphi_{f, 1}\Q$, and have no $1$-eigenspace on $\psi_{f,\neq 1}\Q$ and $\varphi_{f,\neq 1}\Q$.
So (\ref{dist1}) becomes the distinguished triangle
\be\label{dist2} \Q_{D} \to \psi_{f,1}\Q \to \varphi_{f,1}\Q \overset{[1]}{\to} .\ee
The vanishing cycle complex  $\varphi_{f}\Q$ is supported on the singular set $\Sigma$ (which consists of finitely many points), hence so is $\varphi_{f,\neq 1}\Q \simeq \psi_{f,\neq 1}\Q$. Therefore, we have a short exact sequence of MHS  
\be \label{3.9} 
 0 \to H^{n}(D) \to H^{n}(D,\psi_{f,1}\Q) \to H^{n}(D, \varphi_{f,1}\Q) \to 0 
 \ee 
 and MHS isomorphisms 
 \be \label{3.9b} H^{n}(D,\psi_{f,\neq 1}\Q)\simeq \bigoplus _{x\in \Sigma}\mathcal{H}^{n}(\psi_{f,\neq 1}\Q)_{x} \simeq \bigoplus _{x\in \Sigma} H^{n}(F_{x})_{\neq 1},\ee where $H^{n}(F_{x})_{\neq 1}$ denotes the sub-vector space of $H^{n}(F_{x})$ on which the monodromy operator $h^*$ has no $1$-eigenspace.  Note also that for the exactness of (\ref{3.9}) we use the fact that $D$ has only isolated singularities and it is homotopy equivalent to a bouquet of $n$-spheres (as $f$ is transversal at infinity).
 
Therefore the claim about the spectral pairs of $H^n(\U^c)$ in the two cases corresponding to weight $n$ and $n+1$, respectively, follows by combining (\ref{3.1}) and (\ref{3.2}) with (\ref{3.9}) and (\ref{3.9b}). 
\end{proof}


Theorem \ref{t3.2} also shows that $Gr^{W}_{n+1} H^{n}(\U^{c})$ is non-trivial only for eigenvalue $1$, and $Gr^{W}_{n} H^{n}(\U^{c})$ is non-trivial  only for eigenvalues different from $1$. Hence we get:
 
\bc With the same assumptions as in Theorem \ref{t3.2}, we have that \begin{itemize}
\item[(1)] $H^{n}(\U^{c})$ has a pure MHS of weight $n$ if, and only if,  $\delta_{n}(\U,1)\neq 0$.
\item[(2)] $H^{n}(\U^{c})$ has a pure MHS of weight $n+1$ if, and only if,   $\delta_{n}(\U, t)$ is a power of $t-1$.
\end{itemize}
\ec

It was shown in \cite[Proposition 2.1]{Max} and \cite[Corollary 5.4]{Liu} that if $f$ is transversal at infinity and $D=f^{-1}(0)$ is a rational homology manifold, then $\delta_{n}(\U,1)\neq 0$.  Therefore, we obtain the following:
\bc If under the assumptions of Theorem \ref{t3.2} we assume, moreover, that $D$ is a rational homology manifold, then the Alexander module $H^{n}(\U^{c})$ has a pure MHS of weight $n$. 
\ec

\br\label{homog} Let $\widetilde{f}$ be the homogenization of $f$. Then  $\widetilde{f}$ is a homogeneous polynomial with a one-dimensional singular locus in $\C^{n+2}$, whose transversal singularities are in one-to-one correspondence with the singular points of $D$. Here, the one-dimensional singular locus of $\widetilde{f}$ has finitely many irreducible components, and the transversal singularities are defined as the intersection of each irreducible component with a generic hyperplane. Let $\widetilde{F}=\{\widetilde{f}=1\}$ be the (global) Milnor fibre of $\widetilde{f}$ at the origin. As shown by Libgober \cite[Proposition 4.9]{L94} (see also \cite[Corollary 6.5]{Liu}), there is a $\Gamma$-module and MHS isomorphism
$H^{n}(\U^c)\cong H^{n}(\widetilde{F})$, and the monodromy action is canonical. 
So the result in Theorem \ref{t3.2} can be also viewed as describing the spectral pairs for  $H^n(\widetilde{F})$; compare also with \cite[Section 5]{Max}.
\er


\subsection{Plane curves}
In this section, we focus on the case of plane curves, i.e., $n=1$. If $f:\C^2\to \C$ is a reduced degree $d$ polynomial transversal at infinity, then $\overline{D}\cap H_{\infty}$ consists of $d$ distinct points. So, without loss of generality, we can assume that the top degree part of $f$ is given by $f_{d}=x_{1}^{d}+x_{2}^{d}$. Following Steenbrink's formula, an easy computation gives us all the non-zero spectral pairs for $H^1(F_{\infty})$, with $F_{\infty}=\{x_{1}^{d}+x_{2}^{d}=1\}$, namely: 
\begin{itemize} 
\item[(1)] $ h^{0,1}_{\alpha}(F_{\infty})=h^{1,0}_{1-\alpha}(F_{\infty})=d\alpha -1$, for $\alpha=\dfrac{j}{d}$ with $1\leq j \leq d-1$;
\item[(2)] $ h^{1,1}_{0}(F_{\infty})=d-1$.
\end{itemize} 

On the other hand, it is known that $\dim H_{1}(\U^{c})_{1} =r-1$, where $r (\leq d)$ is the number of the irreducible components of $D$, see \cite[Lemma 21]{Oka}. Combing all the above, we obtain the following.
\bc \label{c3.6} Assume that the reduced degree $d$ polynomial $f: \C^{2}\rightarrow \C$ is transversal at infinity. Then we have the following upper bounds for the (possibly) non-zero spectral pairs of $H^1(\U^c)$:
 \begin{itemize}
\item[(1)] $ h^{0,1}_{\alpha}( \U^{c})=h^{1,0}_{1-\alpha}(\U^{c})\leq d\alpha -1$, for $\alpha=\dfrac{j}{d}$ with $1\leq j \leq d-1$;
\item[(2)] $ h^{1,1}_{0}(\U^{c})=r-1 $.
\end{itemize} 
In particular,  $h_{\frac{1}{d}}^{0,1}( \U^{c})=h_{\frac{d-1}{d}}^{1,0}( \U^{c})=0$.
\ec


\subsection{Line arrangements}
In this section, we specialize further to the case when $f$ defines a line arrangement in $\C^2$.

\bd\label{313} For any positive integer $m$, we let $\widehat{m}: (0,1) \to \Z$ be the function defined by:  $$\widehat{m}(\alpha)=\left\{ \begin{array}{ll}
m \alpha , & \text{  if  }  m \alpha \in \Z, \\
1, & \text{ otherwise}.\\
\end{array}\right.  $$
\ed

In \cite{Mas}, Massey obtained bounds on the Betti numbers of the Milnor fibre $F=\{Q=1\}$ of a central hyperplane arrangement in $\C^3$, where $Q$ is a complex homogeneous polynomial of degree $d$ in three variables which can be factored as a product of linear forms. In view of Remark \ref{homog} and with the notations from Definition \ref{313}, the following result can be regarded as a refinement of Massey's result to the level of spectral pairs.

\bt \label{t3.3} Assume that the reduced degree $d$ polynomial $f: \C^{2}\rightarrow \C$ defines a line arrangement $D$, which is transversal at infinity. Say $D$ has $k$ singular points, whose multiplicities are denoted by $m_{1},\cdots,m_{k}$, with $2\leq m_{i} \leq d$. Then we have the following upper bounds for the (possibly) non-zero spectral pairs of $H^1(\U^c)$:
 \begin{itemize}
\item[(1)] $ h^{0,1}_{\alpha}( \U^{c})=h^{1,0}_{1-\alpha}(\U^{c})\leq \min \lbrace d\alpha-1,\sum_{i=1}^{k} (\widehat{ m_{i}}(\alpha)-1) \rbrace$ for $\alpha=\dfrac{j}{d}$ with $1\leq j \leq d-1.$
  
\item[(2)] $ h^{1,1}_{0}(\U^{c})=d-1$.
\end{itemize}  
In particular,  $h_{\frac{j}{d}}^{0,1}( \U^{c})=h_{\frac{d-j}{d}}^{1,0}( \U^{c})=0$ if $\gcd(j,d)=1$, unless $f$ is already homogeneous (up to a change of coordinates).
\et
\begin{proof}
For a singular point $x\in D$ with multiplicity $m$, the corresponding spectral pairs can be computed by Steenbrink's formula for $x_1^{m}+x_2^{m}$. Then the claim follows from Theorem \ref{t3.2}.

If there exists $1\leq j\leq d-1$ such that $\gcd(j,d)=1$ and $h_{\frac{j}{d}}^{0,1}( \U^{c})=h_{\frac{d-j}{d}}^{1,0}( \U^{c})\neq 0$, then  one of $\{ m_{i} \}$ must be $d$. Since $f$ defines a line arrangement of degree $d$, it follows that $f$ is homogeneous, up to a change of coordinates.
\end{proof}

\br
Massey's result can be restated as follows: for any  $\alpha=\dfrac{j}{d}$ with $1\leq j \leq d-1$, we have that 
$$h^{0,1}_{\alpha}(\U^{c})+h^{1,0}_{\alpha}(\U^{c}) =h^{0,1}_{1-\alpha}(\U^{c})+h^{1,0}_{1-\alpha}(\U^{c})\leq \sum (m_{i}-2),$$
where the sum is over all $m_{i}$ such that $ m_{i}\alpha \in \Z$.
\er


\section{Spectral pairs of Alexander modules of boundary manifold}
\subsection{Peripheral complex as a mixed Hodge module}

In \cite{LiM}, we showed that the Alexander modules $H^{i}(M^{c})$ of the boundary manifold $M$ are endowed with mixed Hodge structures induced from the mixed Hodge module structure of the peripheral complex.\footnote{As in the previous section, it is more convenient to work with cohomological invariants of the boundary manifold.} More precisely, in the notations of (\ref{notr}), we have the following:
\bp (\cite[Corollary 3.3]{LiM}) \label{p4.1} Assume that $f:\CN\to \C$ is a reduced degree $d$ polynomial  transversal at infinity.  Then the  peripheral complex $R^{\bullet}$ underlies a (shifted) algebraic mixed Hodge module, hence the Alexander modules of the boundary manifold $M$ carry (non-canonical) MHS.
 \ep
 \begin{proof}
Let us briefly recall the proof. We will use freely the notations of the previous sections. 

It is shown in \cite{Max} that $\R\vert_{H_{\infty} \setminus V\cap H_{\infty}}$ is a local system  $\cL_{\infty}$ with stalk $\Gamma/(t^{d}-1)$ placed in degree $1$. This local system can be obtained from the Milnor fibre $F_{\infty}$, namely: $$\cL_{\infty}\simeq Rp_{\ast} \Q_{F_{\infty}} ,$$ where $p$ is the $d$-fold covering map  $F_{\infty} \to H_{\infty} \setminus V\cap H_{\infty}$. Since $p$ is a finite algebraic map, it follows that $\cL_{\infty}[n]$ is a perverse sheaf on $H_{\infty} \setminus V\cap H_{\infty}$ underlying an algebraic mixed Hodge module.
 Moreover, $$\cL_{\infty}\otimes \C \simeq \bigoplus _{k=0}^{d-1} \cL_{k},$$ where $\cL_{\infty} \otimes \C$ is the complexification of $\cL_{\infty}$, and $\cL_{k}$ is a rank-one $\C$-local system on $H_{\infty} \setminus V\cap H_{\infty}$ with associated representation defined by sending all generators of $H_{1}(H_{\infty} \setminus V \cap H_{\infty}; \Z)$ to $e^{2 \pi i (k/d)}$. In particular, $\cL_{0}$ is the constant sheaf.  Moreover, the associated monodromy action on $\cL_{k}$ has $e^{2 \pi i (k/d)}$ as the only eigenvalue. 
 
 After applying the forgetful functor (\ref{for}), we therefore have the following identity of contructible $\Q$-complexes: 
 \be  \label{4.1} 
 R^{\bullet}\vert_{(H_{\infty}\setminus V\cap H_{\infty})}\simeq \cL_\infty \simeq \Q_{H_{\infty} \setminus V \cap H_{\infty}} \oplus \cL_{\infty, \neq 1}.
  \ee

Let $i_{v}$ and $i_{\infty}$ be the inclusions of $D$ and $H_{\infty}\setminus V\cap H_{\infty}$ into their projective completions $V$ and $H_{\infty}$, respectively. Then the transversality (at infinity) assumption implies that (cf. \cite[Theorem 3.1]{LiM}):
\be \label{4.2a}
R^{\bullet}\vert_{V}\simeq R(i_{v})_{ \ast} \psi_{f}\Q.
\ee
\be  \label{4.2b}
R^{\bullet}\vert_{H_{\infty}}\simeq R(i_{\infty})_{ \ast} \cL_{\infty}.
\ee
Moreover, it was shown in \cite{LiM} that (\ref{4.2a}) and (\ref{4.2b}) underly quasi-isomorphisms of $\Gamma$-complexes. 

By using  of (\ref{4.2a}) and (\ref{4.2b}), there are two ways to proceed in order to show that $R^{\bullet}$ underlies a (shifted) mixed Hodge module. The first possibility  is already described  in the proof of \cite[Corollary 3.3]{LiM} via the quasi-isomorphism (\ref{4.2a}). We present here the second approach, by making use of (\ref{4.2b}).

Consider the inclusions $D  \overset{s}{\hookrightarrow} V \cup H_{\infty} \overset{r}{\hookleftarrow} H_{\infty}$ with $s$ open and $r$ closed,
and the associated distinguished triangle in $D^b_c(V \cup H_{\infty},\Q)$ obtained by using (\ref{pn}) and (\ref{4.2b}):
\be \label{4.3}
s_! \psi_{f}\Q[n]  \to R^{\bullet}[n] \to  Rr_\ast R(i_{\infty})_{ \ast}\cL_{\infty} [n] \overset{[1]}{\to}
\ee
Since $\Q[n+1]$ is a perverse sheaf on $\C^{n+1}$ underlying a mixed Hodge module, and the functor $\psi_f[-1]$ preserves perverse sheaves (and mixed Hodge modules), it follows that $\psi_{f}\Q[n]$ is a perverse sheaf on $D$ underlying  a mixed Hodge module. Since  $s$ is a quasi-finite affine map, it follows from \cite[Corollary 5.2.17]{D2} that $s_!\psi_{f}\Q[n]$ is a perverse sheaf underlying a mixed Hodge module on $V \cup H_{\infty}$.  On the other hand, as $i_\infty$ is a quasi-finite affine morphism, it follows as above that $R(i_{\infty})_{ \ast}\cL_{\infty} [n]$ is a perverse sheaf on $H_{\infty}$ underlying  a mixed Hodge module. Finally, since $r$ is proper, $r_!=r_*$ preserves perverse sheaves and resp. mixed Hodge modules, so $Rr_\ast R(i_{\infty})_{ \ast}\cL_{\infty} [n]$ is a perverse sheaf on $V \cup H_{\infty}$ underlying a mixed Hodge module.
 
So the first and the third terms in the distinguished triangle (\ref{4.3}) are perverse sheaves underlying mixed Hodge modules. In particular, $R^{\bullet}[n]$ can be regarded as an extension of perverse sheaves, both of which underly mixed Hodge modules. We can now argue as in \cite[Corollary 3.3]{LiM} to conclude that there exists a mixed Hodge module $\mathbf{R}$ (defined as an extension of the above mentioned mixed Hodge modules) such that ${\rm For} (\mathbf{R})= R^{\bullet}[n]$, where ${\rm For}:\MHM \to {\rm Perv}_{\Q}$ is the forgetful functor assigning to a mixed Hodge module the corresponding perverse sheaf.

Finally, by Proposition \ref{p2.1}, it follows that the Alexander modules $H^{i}(M^{c})$, for all $i$, are endowed with (non-canonical) mixed Hodge structures induced from (shifted) mixed Hodge module structure of the peripheral complex.
 \end{proof}
 
Consider now the long exact sequence of MHS associated to the distinguished triangle (\ref{4.3}). As shown in  \cite[4-16]{LiM}), one gets a short exact sequence for the middle degree $n$:
\be \label{4.4b}
0 \to H^{n}_{c}(D, \psi_{f}\Q) \rightarrow H^{n}(M^{c}) \rightarrow H^{n}(F_{\infty}) \rightarrow 0.
\ee
Moreover, we have MHS isomorphisms (\cite[4-14,4-15]{LiM})
\be  H^{i}(M^{c}) \cong H^{i}(F_{\infty})  \text{ for } i<n ,
\ee 
and
\be H^{i}(M^{c}) \cong H^{i}_{c}(D, \psi_{f}\Q)  \text{ for } i>n .
\ee
Since the nearby cycle complex  $\psi_f \Q$ is (up to a shift) self-dual as a mixed Hodge module,   
the MHS on $ H^{i}_{c}(D, \psi_{f}\Q) $ and $H^{2n-i}(D, \psi_{f}\Q) $ are dual to each other at level $n$ (in the sense of \cite[Definition 1.6.1]{F}).

Combing the above results with the MHS isomorphism (\ref{3.0}), one gets that the MHS on $H^i(M^c)$
and $H^{2n-i}(M^c)$ are dual to each other at level $n$ for $i\neq n$, and there exist MHS isomorphisms $$H^i(M^c)\cong H^i(\U^c)$$ for $i<n$, which can be viewed as being induced by the $n$-homotopy equivalence $M\hookrightarrow \U$. This shows that there exists a unique choice for the MHS on $H^{i}(M^{c})$ for $i\neq n$. Moreover, this mixed Hodge structure on $H^{i}(M^{c})$ ($i\neq n$) is compatible with the monodromy action. Note that Theorem \ref{LiM} asserts that $H^{i}(M^{c})$ ($i\neq n$) is in fact semi-simple.

\br  Here, we define the MHS on $H^i(M^c)$ as follows: 
$$H^i(M^c) \cong H ^i(V\cup H_\infty, \mathbf{R}), $$
such that the MHS on $H^i(M^c)$ is isomorphic with the one on $H^i(\U^c)$ for $i<n$.
The corresponding MHS for the homology version is $H_i(M^c)\cong H ^i(V\cup H_\infty, \mathbf{R})(n) $, where $(n)$ denotes the Tate twist by $n$.
\er

\br \label{r4.2} Since the mixed Hodge module $\mathbf{R}$ inducing a MHS on the Alexander modules of the bounday manifold $M$ is defined as an extension of mixed Hodge modules, it is not clear to us at this point if there exists a unique choice of MHS on $H^{n}(M^{c})$.

In fact, as already mentioned in the proof of Proposition \ref{p4.1}, there are two ways to construct the underlying mixed Hodge module for peripheral complex, by making use respectively of either one of the quasi-isomorphisms (\ref{4.2a}) and (\ref{4.2b}). In our previous work \cite{LiM}, we have used (\ref{4.2a}) to do this and in this case we have the following short exact sequence of MHS: \be \label{4.4a} 
0 \to H^{n}_{c}(F_{\infty}) \rightarrow H^{n}(M^{c}) \rightarrow H^{n}(D, \psi_{f}\Q) \rightarrow 0.
\ee  
The short exact sequences (\ref{4.4b}) and (\ref{4.4a}) have dual MHS at level $n$ on the first term and third terms.  Since the middle terms in these two short exact sequences coincide, the mixed Hodge numbers associated to $H^{n}(M^{c})$ are well-defined  if and only if the MHS $H^{n}(M^{c})$ is self-dual. This self-duality of the  MHS $H^{n}(M^{c})$ would follow if one can show that the underlying mixed Hodge module of the peripheral complex is self-dual as a mixed Hodge module. While this is still an open question in general, we manage to give an affirmative answer in the case of hypersurfaces which are either plane curves (e.g., line arrangements) or rational homology manifolds with only isolated singularities, see  Propositions \ref{p4.7} and \ref{p4.8} below.
\er


\subsection{Isolated singularities} In this section, we assume that the hypersurface $D=\{f=0\}$ has only isolated singularities. 
We recall the following result about the (only non-trivial) Alexander polynomial of $M$:
 \bp\label{p4.3} (\cite[Corollary 1.3]{LiM}) Assume that $f:\CN \to \C$ is a reduced degree $d$ polynomial transversal at infinity, and $D=f^{-1}(0)$ has at most isolated singularities. Let $\Sigma$ denote the singular locus of $D$.  Then the only non-trivial Alexander module of $M$ appears in degree $n$, and we have the following polynomial identity:
\begin{equation}
(t-1)^{(-1)^{n+1}+\mu}(t^{d}-1)^{\xi}\cdot \prod_{x\in \Sigma}\Delta_{x}(t)=\delta_{n}(M,t)= \delta_{n}(\U, t)^{2} \cdot e(t)
\end{equation}
 where $ \Delta_{x}(t)$ is the top local Alexander polynomial associated to the point $x\in \Sigma$, $\xi=\dfrac{(d-1)^{n+1}+(-1)^{n}}{d}$, and  $\mu=(d-1)^{n+1} - \sum_{x\in \Sigma} \mu_{x}$. Here $\mu_{x}$ is the Milnor numer of $f$ at $x$.  Moreover,  the degree of $e(t)$ is always even and the degree of $\delta_{n}(M,t)$ is  $2(d-1)^{n+1}$, which does not depend on the type (nor position) of singularities. 
\ep
\br \label{r4.4}  $e(t)$ is called the error term since it is the product of the error terms in the divisibility results (\ref{d1}) and (\ref{d2}). 
\er

Since the peripheral complex $\cR^{\bullet}$ is a torsion $\Gamma$-module,  the monodromy $t$-action give rise (in the notation of (\ref{notr})) to a decomposition:
 \be R^{\bullet}\simeq R^{\bullet}_{1} \oplus R^{\bullet}_{\neq 1}.\ee
In what follows, we investigate separately the unipotent (eigenvalue-$1$) part $R^{\bullet}_{1}$ and, respectively, the non-unipotent part $R^{\bullet}_{\neq 1}$ corresponding to eigenvalues other than $1$.


\subsubsection{Non-unipotent part}
 Note that $\psi_{f,\neq 1}\Q \simeq \varphi_{f,\neq 1}\Q$ is supported on the finite set of points  $\Sigma=Sing(D)$. Then  
 \be\label{nj} R(i_{v})_{ \ast} \psi_{f, \neq 1}\Q \simeq R(i_{v})_{  \ast} \varphi_{f, \neq 1}\Q \simeq (i_{v})_{  !} \varphi_{f, \neq 1}\Q\ee is also supported on $\Sigma$. For the second quasi-isomorphism in (\ref{nj}), if $k:\Sigma \hookrightarrow D$ denotes the (closed) inclusion map, the support condition for $\varphi_{f,\neq 1}\Q$ yields that $\varphi_{f,\neq 1}\Q\simeq Rk_*k^*\varphi_{f,\neq 1}\Q\simeq k_!k^*\varphi_{f,\neq 1}\Q$. Since $l:=i_v \circ k:\Sigma \hookrightarrow V$ is also closed, we then get: $R(i_{v})_{  \ast} \varphi_{f, \neq 1}\Q \simeq Rl_*k^*\varphi_{f, \neq 1} \simeq l_! k^*\varphi_{f, \neq 1} \simeq (i_v)_! k_!k^*\varphi_{f, \neq 1} \simeq (i_{v})_{  !} \varphi_{f, \neq 1}$.

 The quasi-isomorphisms (\ref{4.1}) and (\ref{4.2a}) (which underly quasi-isomorphisms of $\Gamma$-complexes) give us that:
  $$ R^{\bullet}_{\neq 1} \vert_{ V} \simeq  R(i_{v})_{  \ast} \psi_{f, \neq 1}\Q \ \text{   and  } \ R^{\bullet}_{\neq 1}\vert_{(H_{\infty}\setminus V\cap H_{\infty})}\simeq  \cL_{\infty, \neq 1}.$$ Note that $\cL_{\infty, \neq 1}$ and $R(i_{v})_{  \ast} \psi_{f, \neq 1}\Q$ have disjoint supports, hence if $a:H_{\infty} \setminus V \cap H_{\infty} \hookrightarrow V \cup H_{\infty}$ and $b :V \hookrightarrow V \cup H_{\infty}$ denote the inclusion maps, the associated distinguished triangle
  $$a_!\cL_{\infty, \neq 1} \to R^{\bullet}_{\neq 1} \to b_*R(i_{v})_{  \ast} \psi_{f, \neq 1}\Q \overset{[1]}{\to}$$ splits (e.g., see \cite[Corollary 1.2.7]{Ne}), i.e., 
  \be \label{4.6}
  R^{\bullet}_{\neq 1}  \simeq \cL_{\infty, \neq 1} \oplus \varphi_{f,\neq 1} \Q ,
  \ee
 where  $\cL_{\infty, \neq 1} $ and $ \varphi_{f,\neq 1} \Q $ are viewed as complexes of sheaves on $V\cup H_{\infty}$ after applying the extension by $0$ from their supports.   In particular, $R^{\bullet}_{\neq 1}$ underlies (up to a shift) a unique self-dual mixed Hodge module.
  
We can now prove the following: 
\bt \label{t4.5}  Assume that $f:\CN\to \C$ is a reduced degree $d$ polynomial transversal at infinity so that $D=f^{-1}(0)$ has at most isolated singularities.  Let $\Sigma$ denote the singular locus of $D$. Then  we have  MHS and $\Gamma$-module isomorphisms:
\be  \label{4.7}
 H^{n}(M^{c})_{\neq 1} \cong H^{n}(F_{\infty})_{\neq 1} \oplus \big(  \bigoplus_{x\in \Sigma} H^{n}(F_{x})_{\neq 1} \big), 
 \ee
 with $F_x$ denoting the Milnor fiber of $f$ at $x \in \Sigma$.  Therefore, for $\alpha>0$, we get:
 \be 
 h^{p,q}_{\alpha}(M^{c}):=h^{p,q}_{\alpha}(H^{n}(M^{c}))=\sum_{x\in \Sigma} h_{\alpha}^{p,q}(F_{x}) + h^{p,q}_{\alpha}(F_{\infty})  .
 \ee   In particular, $ H^{n}(M^{c})_{\neq 1}$ with the induced $t$-action is a symmetric MHS at level $n$, hence it is self-dual at level $n$. 
\et  
 \begin{proof} 
 It is clear that the quasi-isomorphism (\ref{4.6}) yields $$ H^{n}(M^{c})_{\neq 1} \cong H_c^{n}(F_{\infty})_{\neq 1} \oplus \big(  \bigoplus_{x\in \Sigma} H^{n}(F_{x})_{\neq 1} \big) .$$
 Then (\ref{4.7}) follows from the duality isomorphism  (of MHS and semi-simple $\Gamma$-modules) $H^{n}_{c}(F_{\infty})_{\neq 1} \cong H^{n}(F_{\infty})_{\neq 1} $.  
 
 The MHS $ H^{n}(M^{c})_{\neq 1}$  is symmetric at level $n$ since the same is true for the summands on the right-hand side of (\ref{4.7}), corresponding to the Milnor fibre $F_{x}$ of $f$ at $x \in \Sigma$ and, resp., the Milnor fibre at infinity $F_{\infty}$.
 \end{proof}

\br
The monodromy theorem tells us that the size of Jordan blocks for the monodromy action on $H^{n}(F_{x})_{\neq 1}  $ is at most $n+1$, hence the same is true for $H^{n}(M^{c})_{\neq 1}$. 
\er 


\subsubsection{Unipotent part}
Unlike the non-unipotent part, the corresponding MHS on $H^{n}(M^{c})_{1}$ with induced monodromy action is in general not symmetric for weight $n$, see Remark \ref{r4.9} below. In particular, the relation between the relative weight filtration and the nilpotent operator $N$ is in general more complicated. Nevertheless, for hypersurfaces in general position at infinity which are, moreover, rational homology manifolds, we obtain a complete characterization of the mixed Hodge numbers of $H^{n}(M^{c})_{1}$, as shall be explained below.

  Let $v$ denote the inclusion $V \hookrightarrow \CP$. Then the induced map $v^k: H^k(\CP)\to H^k(V)$ is a monomorphism for all $k$ with $0\leq k \leq 2n$ (e.g., see \cite[(5.2.17)]{D1}). The {\it primitive $k$-th cohomology group} of $V$ is defined as $$H^k_0(V):= \mathrm{coker} (v^k).$$ In particular, for any $0\leq k \leq 2n$,  we have:
  \be  \label{primitive}
   h^{p,q}(H^k_0(V))= h^{p,q}(H^k(V))-h^{p,q}(H^k(\CP)).
  \ee

\bp \label{p4.7} Assume that $f:\CN\to \C$ is a reduced degree $d$ polynomial transversal at infinity so that $D=f^{-1}(0)$ has at most isolated singularities. Let $\Sigma$ denote the singular locus of $D$.  If $D$ is a rational homology manifold (e.g., $D$ is smooth), then $H^{n}(M^{c})_{1}$ has only $3$ possible weights: $(n-1)$, $n$ and $(n+1)$, and we have:
 $$h^{p,q}_{0}(H^{n}(M^{c}))=
\begin{cases}
h^{p,q}_0 (F_{\infty}),           & p+q=n+1,  \\
h^{p,q}(H^{n}_{0}(V_{sm})) - \sum_{x\in \Sigma} \dim Gr_{p}^{F} H^{n}(F_{x}), & p+q=n, \\
h^{p,q}(H_{0}^{n-1}(V\cap H_{\infty})), & p+q=n-1, \\
\end{cases} 
$$
where $H^{n-1}_{0}(V\cap H_{\infty})$ denotes the corresponding primitive cohomology for the hypersurface $V\cap H_\infty$ in $H_\infty \cong \mathbb{CP}^n$, and $V_{sm}$ is any smooth degree $d$ hypersurface in $\CP$. 
Moreover,  $$h^{p,n-1-p}(H_{0}^{n-1}(V\cap H_{\infty}))=h^{p+1,n-p}_0 (F_{\infty})= h^{n-p, p+1}_0(F_{\infty}),$$ hence the mixed Hodge numbers associated to $H^{n}(M^{c})$ are self-dual.
\ep

\begin{proof}
 If $D$ is a rational homology manifold, then $\varphi_{f,1} \Q=0$, hence by (\ref{dist2}) we get \be\label{n33} H^{n}_{c}(D, \psi_{f,1}\Q) \cong H^{n}_{c}(D).\ee
Moreover, the Poincar\'e duality isomorphism holds for $D$. Since $D\simeq \bigvee_\mu S^{n}$, we then have $$H^{i}_{c}(D)=\left\{ \begin{array}{ll}
\mu  , & i=n, \\
1 , & i=2n \\
0, & \text{ otherwise}.\\
\end{array}\right.$$ 
 
Consider now the following long exact sequence of MHS for $n>1$:
\be\label{ref} 0 \to H^{n-1}(V)\to  H^{n-1}(V\cap H_{\infty})\to H^{n}_{c}(D) \rightarrow H^{n}(V) \to H^{n}(V\cap H_{\infty})\to 0.\ee
If $n=1$, this long exact sequence of MHS still exists, since the transversality assumption implies that  $V\cap H_\infty$  is a finite set, hence the last term $H^{1}(V\cap H_{\infty})=0$.

The transversality assumption implies that $V$ is also a rational homology manifold, hence $H^{n}(V)$ has a pure Hodge structure of weight $n$ (see \cite{De}). Moreover, the corresponding Hodge numbers $h^{p,n-p}(H^{n}(V))$ can be computed from Hodge-theoretic invariants associated to the local Milnor fibers $F_{x}$, for $x\in \Sigma$, as follows (see \cite{D96}): \be\label{hn} h^{p,n-p}(H^{n}(V))=h^{p,n-p}(H^{n}(V_{sm}))- \sum_{x\in \Sigma} \dim Gr_{p}^{F} H^{n}(F_{x}) ,\ee
where $V_{sm}$ is a smooth degree $d$ hypersurface in $\CP$ and its Hodge numbers are computable only in terms of $d$ and $n$. Note that, by transversality, $V\cap H_{\infty}$ is smooth, hence $H^{n-1}(V\cap H_{\infty})$ has a pure Hodge structure of weight $n-1$. 
 So we see by (\ref{ref}) that $H^{n}_{c}(D)$ is a MHS with possible weights $n$ and $n-1$.
 
Since $V$ has only isolated singularities and $V\cap H_{\infty}$ is smooth by transversality, we have by \cite[(5.2.6),(5.2.11)]{D1} the isomorphisms: $ H^{n-1}(V) \cong H^{n-1}(\CP)\cong H^{n-1}(\mathbb{CP}^{n})$ and $H^{n}(V\cap H_{\infty}) \cong H^{n}(\mathbb{CP}^{n}) \cong  H^{n}(\mathbb{CP}^{n+1})$. So the contribution of $H^{n-1}(V\cap H_{\infty})$ to the weight $n-1$ part of $H_{c}^{n}(D)$ comes in the form of the primitive cohomology group $H_{0}^{n-1}(V\cap H_{\infty})$. On the other hand, the contribution of $H^n(V)$ to $h^{p,q}(H_{c}^{n}(D))$ with $p+q=n$ comes in the form of $h^{p,q}(H^n(V))-h^{p,q}(H^{n}(V\cap H_{\infty}))$, which by the  above isomorphisms and (\ref{primitive}) equals $h^{p,q}(H_0^n(V))$. Finally, $h^{p,q}(H_0^n(V))$ can be computed from (\ref{hn}) by subtracting $h^{p,q}(H^n(\CP))$ on both sides and using (\ref{primitive}). Altogether, we get that
     $$h^{p,q}(H_{c}^{n}(D))=\left\{ \begin{array}{ll}
h^{p,n-p}(H^{n}_{0}(V_{sm})) - \sum_{x\in \Sigma} \dim Gr_{p}^{F} H^{n}(F_{x}), & p+q=n, \\
h^{p,q}(H_{0}^{n-1}(V\cap H_{\infty})), & p+q=n-1, \\
0, & \text{otherwise,}\\
\end{array}\right.$$

The claim about $h^{p,q}_{0}(H^{n}(M^{c}))$ follows now from (\ref{n33}) and the additivity property for the spectral pairs associated to the unipotent part of the short exact sequence (\ref{4.4b}), after recalling that $H^n(F_{\infty})$ has only weights $n$ and $n+1$, with the weight $n+1$ corresponding exactly to the unipotent part $H^n(F_{\infty})_1$, see (\ref{36}).
 
 Finally, as shown in \cite[Section 4]{Ste2}, one has a MHS isomorphism 
 $$Gr^W_{n+1} H^{n}(F_{\infty}) \cong H^{n-1}_{0}(V\cap H_{\infty})(-1) .$$ 
So the equality $h^{p,n-1-p}(H_{0}^{n-1}(V\cap H_{\infty}))=h_0^{p+1,n-p} (F_{\infty})= h_0^{n-p, p+1}(F_{\infty})$ follows since the weight $n+1$ of $H^{n}(F_{\infty})$ is given exactly by the unipotent part $H^n(F_{\infty})_1$.  So the MHS on $H^{n}(M^{c})_{1}$ is self-dual at level $n$. Since this is also true for $H^{n}(M^{c})_{\neq 1}$ by Theorem \ref{t4.5}, the claim that $H^n(M^c)$ is self-dual follows. 
\end{proof}

\br   The Hodge numbers $ h^{p,q}(H^{n}_0(V_{sm}))$ of the above result can also be computed by using the Milnor algebra $\M(n,d)$. In fact,  there exists an MHS isomorphism (see \cite[Section 4]{Ste2})
 $$Gr^W_n(H^n(F_\infty)) \cong H^n_0 (V_{sm}) ,$$
 where we choose $V_{sm}$ to be the projective closure of the Milnor fibre $F_\infty$, which is smooth. Therefore,
 $$h^{p,n-p}(H^{n}_0(V_{sm}))= \sum_{i=1}^{d-1} \dim \M(n,d)_{pd+i-n-1}- h^{p,n-p}(H^n(\CP)).$$
\er


\subsection{Plane curves}
In this section, we focus on the case of plane curves (i.e., $n=1$). 
We will use freely the notations and assumptions of the previous sections. In particular, $f:\C^2 \to \C$ is a degree $d$ reduced polynomial in general position at infinity, with $D=f^{-1}(0)$, $\Sigma=Sing(D)$, and $V$ is the projective completion of $D$.

Assume that $D$ has $r$ irreducible components, and the germ $D_{x}$ at a singular point $x$ has $r_{x}$ irreducible local branches. Note that $r_{x}$ may be different from the number of irreducible components of $D$ passing through $x$, denoted by $n_{x}$. 

The MHS of $H^{\ast}(V)$ are determined by the local singularities, the number of irreducible components of $V$ and the degree of $f$. An explicit formula for the corresponding mixed Hodge numbers is given in \cite[Proposition 2.2]{Abd} as follows:
\begin{itemize}
\item[(1)] $h^{0,0}(H^{0}(V))=1$.
\item[(2)]$h^{0,0}(H^{1}(V))=  \sum_{x\in \Sigma} (r_{x}-1)+1-r$.
\item[(3)]$h^{0,1}(H^{1}(V))=h^{1,0}(H^{1}(V))=\dfrac{1}{2}\{\mu +2r-d-1 -\sum_{x\in \Sigma} (r_{x}-1)\}$. 
\item[(4)] $h^{1,1}(H^{2}(V))= r$.
\end{itemize}  

\br The above formulae for $h^{0,0}(H^{1}(V))$ and $h^{0,1}(H^{1}(V))$ look slightly different from the ones in \cite[Proposition 2.2]{Abd}.  In fact, the formula in loc.cit. for $h^{0,0}(H^{1}(V))$ is: \be\label{abd} h^{0,0}(H^{1}(V))= \sum _{x\in \Sigma} \sum_{1\leq j \leq r} (r_{j,x}-1)+  \sum_{x\in \Sigma} (n_{x}-1)+1-r,\ee
where $r_{j,x}$ denotes the number of local branches from $D_{j}$ passing through $x$. It is easy to check that (\ref{abd}) coincides with our formula above. Also, our formula for $h^{0,1}(H^{1}(V))$ follows from loc.cit. together with the fact that $\chi(V)=\chi(D)+d= d+1-\mu$ and the well-known formula for $\chi(V)$ (e.g., see \cite[(5.4.4)]{D1}). (Recall also that by the transversality assumptions, $D$ has the homotopy type of a wedge of $\mu$ circles.)
\er 

By the transversality assumption, $V \cap H_{\infty}$ consists of $d$ distinct points. Then one gets all the non-trivial mixed Hodge numbers of $H^{\ast}_{c}(D)$ by the corresponding  compactly supported long exact sequence, namely: 
\begin{itemize}
\item[(1)]$h^{0,0}(H^{1}_{c}(D))=  \sum_{x\in \Sigma} (r_{x}-1)+d-r$.
\item[(2)]$h^{0,1}(H^{1}_{c}(D))=h^{1,0}(H^{1}_{c}(D))=\dfrac{1}{2}\{\mu +2r-d-1 -\sum_{x\in \Sigma} (r_{x}-1)\}$. 
\item[(3)] $h^{1,1}(H^{2}_{c}(D))= r$.
\end{itemize}

Let us next consider the compactly supported hypercohomology long exact sequence associated to the distinguished triangle (\ref{dist2}):
\be \label{4.10}
0 \to H^{1}_{c}(D) \to H^{1}_{c}(D, \psi_{f,1}\Q) \to H^{1}(\Sigma, \varphi_{f,1}\Q) \to H^{2}_{c}(D) \to H^{2}_{c}(D, \psi_{f,1}\Q)  \to 0. 
\ee
Note that the $\Q$-vector spaces  $H^{1}(\Sigma, \varphi_{f,1}\Q)$, $ H^{2}_{c}(D) $ and $H^{2}_{c}(D, \psi_{f,1}\Q) $ have pure Hodge structures of type $(1,1)$, and of corresponding dimensions $\sum_{x\in \Sigma}(r_{x}-1)$, $r$ and $1$, respectively. In fact, we have 
$$H^{1}(\Sigma, \varphi_{f,1}\Q)\cong \bigoplus_{x\in \Sigma}  H^{1}(F_{x})_{1} ,$$
 and the dimension calculation $\dim H^{1}(F_{x})_{1} =r_{x}-1$ follows from \cite[Theorem 4.2]{DN} and the fact that the Jordan blocks for $H^{1}(F_{x})$ corresponding to the eigenvalue $1$ must all have size $1$.  Moreover, we have by \cite[page 404]{AGV} that $H^{1}(F_{x})_1$ has a pure Hodge structure of type $(1,1)$. By duality, we also have MHS isomorphisms  
$$ H^{2}_{c}(D, \psi_{f,1}\Q) \cong H^{0}(D, \psi_{f,1}\Q) (-1)\cong H^{0}(D)(-1),$$ 
where $(-1)$ denotes the Tate twist (and the second isomorphism follows from (\ref{dist2}) and  the connectivity of the local Milnor fibers). Then one can compute the mixed Hodge numbers associated to $H^{1}_{c}(D, \psi_{f,1}\Q)$ by using (\ref{4.10}).

Consider now the unipotent (eigenvalue-$1$) part of  the  short exact sequence (\ref{4.4b}):
\be  \label{4.9}
0 \to H^{1}_{c}(D, \psi_{f,1}\Q) \rightarrow H^{1}(M^{c})_{1} \rightarrow  H^{1}(F_{\infty})_{1} \rightarrow 0,
\ee
where $H^{1}(F_{\infty})_{1}$ has a pure Hodge structure of type $(1,1)$ and dimension $d-1$.
By combining the above facts with Theorem \ref{t4.5}, we get the following: 

\bp \label{p4.8} Assume that the reduced degree $d$ polynomial $f: \C^{2}\rightarrow \C$ is transversal at infinity. Then the non-zero spectral pairs of $H^{1}(M^{c})$ are listed as follows:\begin{itemize}
 \item[(1)] For $\alpha=0$, $$h^{0,0}_{0}(H^{1}(M^{c}))=h^{1,1}_{0}(H^{1}(M^{c}))= \sum_{x\in \Sigma} (r_{x}-1)+d-r,$$
$$h^{0,1}_{0}(H^{1}(M^{c}))=h^{1,0}_{0}(H^{1}(M^{c})=\dfrac{1}{2}\{\mu +2r-d-1 -\sum_{x\in \Sigma} (r_{x}-1)\}.$$
\item[(2)] For $\alpha>0$, $$h^{0,1}_{\alpha}(H^{1}(M^{c}))=h^{1,0}_{1-\alpha}(H^{1}(M^{c}))=\sum_{x\in \Sigma} h_{\alpha}^{0,1}(F_{x}) + \widehat{d}( \alpha)-1.$$
$$h^{0,0}_{\alpha}(H^{1}(M^{c}))=h^{1,1}_{\alpha}(H^{1}(M^{c}))=\sum_{x\in \Sigma} h_{\alpha}^{0,0}(F_{x}) .$$
\end{itemize}
 In particular, the spectral pairs of $H^{1}(M^{c})$  are determined by the local singularities. Moreover, the mixed Hodge numbers associated to $H^{1}(M^{c})$ are self-dual at level $1$.
\ep  

\br \label{r4.9} Note that $b_{1}(M)=b_{1}(\U)+b_{2}(\U)= 2r+\mu -1$ (\cite[Corollary 2.6]{CoS1}). By \cite[Theorem 4.2]{DN}, $$J_{1}(H^{1}(M^{c}))= b_{1}(M)-1=2r-2+\mu,$$  where $J_1(H^{1}(M^{c}))$ is the number of Jordan blocks associated to $H^{1}(M^{c})$ with eigenvalue 1. One can check that if $H^{1}(M^{c})_{1}$ with the induced $t$-action is a symmetric MHS at level $1$, then $r=1$, i.e., $D$ is irreducible.
\er


\subsection{Line arrangements}
For a line arrangement, let us refer to the following as
the {\it weak combinatorial data}: $d=$ the number of lines in $D$, and $m_{i}=$ the
numbers of points of multiplicity $i$ in $D$, for all $ i \geq 2$. 
In the line arrangement case, $\mu= \sum_{i=1}^{k} (m_{i}-1) -(d-1)$. Hence $H^{1}_{c}(D)$ has a pure Hodge structure $(0,0)$.


\bc \label{c4.11} Assume that the reduced degree $d$ polynomial $f: \C^{2}\rightarrow \C$ defines a line arrangement $D$, which is transversal at infinity, and so that $D$ has $k$ singular points with multiplicities $m_{1},\cdots,m_{k}$, where $2\leq m_{i} \leq d$.  Then the non-zero spectral pairs of $H^{1}(M^{c})$ are determined by the weak combinatorial data and  can be listed as follows:
 \begin{itemize}
 \item[(1)]$ h^{0,0}_{0}(M^{c})=h^{1,1}_{0}(M^{c})=\sum_{i=1}^{k}(m_{i}-1)$.
\item[(2)] $ h^{0,1}_{\alpha}( M^{c})=h^{1,0}_{1-\alpha}(\U^{c})=  \sum_{i=1}^{k} (\widehat{ m_{i}}(\alpha)-1) + \widehat{d}(\alpha)-1$, \  {\rm  for}  $\alpha >0$.
\end{itemize}  
\ec


\end{document}